\documentclass[10pt]{article}
\usepackage{amsfonts}
\usepackage{amsfonts}
\usepackage{amsfonts}
\usepackage{amsmath}
\usepackage{amssymb}
\usepackage{amsfonts}
\usepackage{amsthm}
\usepackage{mathrsfs}
\usepackage[all]{xy}
\usepackage{graphicx}
\begin{document}

\title{Littlewood-Paley characterization for $Q_{\alpha}(\mathbb{R}^n)$ spaces}         
\author{Qifan Li}        
\date{March, 2009}          
\maketitle


\vskip 0.3 cm


\vskip 3 cm

 {\bf Abstract.} In Baraka's paper [2], he obtained the
Littlewood-Paley characterization of Campanato spaces
$L^{2,\lambda}$ and introduced $\mathcal {L}^{p,\lambda,s}$ spaces.
He showed that $\mathcal {L}^{2,\lambda,s}=(-\triangle
)^{-\frac{s}{2}}L^{2,\lambda}$ for $0\leq\lambda<n+2$. In [7], by
using the properties of fractional Carleson measures, J Xiao proved
that for $n\geq2$, $0<\alpha<1$. $(-\triangle
)^{-\frac{\alpha}{2}}L^{2,n-2\alpha}$ is essential the
$Q_{\alpha}(\mathbb{R}^n)$ spaces which were introduced in [4]. Then
we could conclude that $Q_{\alpha}(\mathbb{R}^n)=\mathcal
{L}^{2,n-2\alpha,\alpha}$ for $0<\alpha<1$. In fact, this result
could be also obtained directly by using the method in [2]. In this
paper, We proved this result in the spirit of [2]. This paper could
be considered as the supplement of Baraka's work [2].

\vskip 3 cm

\section{Introduction}

The $Q_{\alpha}$ spaces were first introduced in [1] as a proper
subspace of BMOA defined by means of modified Garcia norm. In [5],
authors showed that: Let $\alpha\in(0,1)$, an analytic function $f$
in the Hardy space $H^1$ on the unit disc belongs to $Q_{\alpha}$,
if and only if its boundary values on the unit circle $\mathbb{T}$
satisfies:

$$\sup_I|I|^{-\alpha}\int_I\int_I\frac{|f(e^{i\theta})-f(e^{i\varphi})|^2}{|e^{i\theta}-e^{i\varphi}|^{2-p}}d\theta
d\varphi<\infty$$ Where the supremum is taken over all subarcs
$I\subset\mathbb{T}$. In [4], the $Q_{\alpha}$ was extended to
Euclidean space $R^n$($n\geq2$). They gave the definition of this
kind of space as follows: For $\alpha\in(-\infty,+\infty)$, $f\in
Q_{\alpha}(\mathbb{R}^n)$ if and only if

$$
\|f\|_{Q_{\alpha}}\triangleq[\sup_Il(I)^{2\alpha-n}\int_I\int_I\frac{|f(x)-f(y)|^2}{|x-y|^{2\alpha+n}}dx
dy]^{\frac{1}{2}}<\infty. \eqno(1.1)$$

Here $I\subset\mathbb{R}^n$ be a cube with the edge parallel to the
coordinate axes, and let $l(I)$ be the length of $I$. The supremum
is taken over all cubs $I\subset\mathbb{R}^n$. There are systematic
research of $Q_{\alpha}(\mathbb{R}^n)$ in [4].

In [4], we have known that if $\alpha<0$, $Q_{\alpha}=BMO$. And if
$\alpha\geq1$, $Q_{\alpha}=\{\mathrm{constants}\}$. We have also
known ([7] theorem 1.2 (1))
$$Q_{\alpha}(\mathbb{R}^n)=(-\triangle
)^{-\frac{\alpha}{2}}L^{2,n-2\alpha}$$for the nontrivial case
$\alpha\in(0,1)$. $L^{2,n-2\alpha}$ denote the Campanato
spaces:$$L^{2,n-2\alpha}\triangleq(\sup_Il(I)^{2\alpha-n}\int_I|f(x)-f_I|^2dx)^{\frac{1}{2}}<\infty.
$$ Combining this result with ([2], theorem 10). We can immediately
obtain:$$Q_{\alpha}(\mathbb{R}^n)=\mathcal
{L}^{2,n-2\alpha,\alpha}$$The Littlewood-Paley characterization is
now clear by the $ \mathcal {L}^{2,n-2\alpha,\alpha}$'s definition
([2], definition 2):
$$\|f\|_{\mathcal {L}^{2,n-2\alpha,\alpha}}\triangleq\sup_I(\frac{1}{|I|^{1-\frac{2\alpha}{n}}}\sum_{j\geq-\log_2l(I)}2^{2\alpha
j}\|\Delta_jf\|_{L^2(I)}^2)^{\frac{1}{2}}\eqno(1.2)$$In this paper
we present an alternative proof of the result. Unlike J. Xiao's
arguments, which make a systematic research of fractional Carleson
measures [3]. Our methods are in the spirit of [2]. We directly
prove the Littlewood-Paley characterization from (1.1) which is the
definition of $Q_{\alpha}(\mathbb{R}^n)$.

Let $\psi(x)$ be a Schwartz function. supp$
\hat\psi(\xi)=\{\xi\in\mathbb{R}^n:\frac{1}{2}\leq|\xi|\leq2\}$ is
compact and $\sum \hat\psi_j(\xi)\equiv1$. We define the
Littlewood-Paley operator by
$$\Delta_j(f)(x)=\psi_j*f(x)$$where$$\psi_j(x)=2^{jn}\psi(2^jx)$$
In this paper, we study the case $f\in S^\prime/\mathcal {P}$. The
homogeneous decomposition of $f$ is given by the formula
$$f=\sum_{j\in \mathbb{Z}}\Delta_j(f)(x)$$
We denote $A\lesssim B$ if $A\leq C(n,\alpha)B$. And define
$A\thickapprox B$ if $A\leq C(n,\alpha)B$ and $B\leq C(n,\alpha)A$.
We have the following main result.

\vskip 0.2 cm \noindent {\bf Main Theorem}  Let $f\in
L^2(\mathbb{R}^n)$, $0<\alpha<1$. We have the Littlewood-Paley
characterization of
$Q_{\alpha}(\mathbb{R}^n)$:$$\|f\|_{Q_{\alpha}}\approx\sup_I(\frac{1}{|I|^{1-\frac{2\alpha}{n}}}\sum_{j\geq-\log_2l(I)}2^{2\alpha
j}\|\Delta_jf\|_{L^2(I)}^2)^{\frac{1}{2}}\eqno(1.3)$$ The main
theorem essentially contains two statements as follows:

\vskip 0.05 cm
If $f\in\mathcal {L}^{2,n-2\alpha,\alpha}$ then
$\|f\|_{Q_{\alpha}}\lesssim\|f\|_{\mathcal
{L}^{2,n-2\alpha,\alpha}}$;

If $f\in Q_{\alpha}(\mathbb{R}^n)$, then $\|f\|_{\mathcal
{L}^{2,n-2\alpha,\alpha}}\lesssim\|f\|_{Q_{\alpha}}$. \vskip 0.05 cm
\vskip 0.15 cm
 \noindent {\bf Remark }  From the main theorem, we get
the relationship between $Q_{\alpha}$ spaces and Morrey type Besov
spaces: In [6], authors introduced a kind of Morrey type Besov
spaces:$$\|f\|_{MB^{p,\sigma}_{\alpha,q}}\triangleq(\sum_{j\in\mathbb{Z}}(\sup_I\frac{1}{|I|^{\frac{\sigma}{n}}}\int_I(2^{\alpha
j}|\Delta_jf|)^qdx)^{\frac{p}{q}})^{\frac{1}{p}}<\infty$$ We
immediately have the embedding property:
$MB^{2,n-2\alpha}_{\alpha,2}\subset Q_{\alpha}$ for $0<\alpha<1$.

\section{Preliminary Lemmas}

The proof of the main theorem relies on following lemmas. To start
with, we introduce some notations: Let $I$ be the any fixed cub in
$\mathbb{R}^n$ with the edge parallel to the coordinate axes. We let
$D_k(I)$, $k\geq0$, denote the set of the $2^{kn}$ subcubes of edge
length $2^{-k}l(I)$ obtained by $k$ successive bipartition of each
edge of $I$. We define $D(I)$ be the set of all the dyadic subcubes
of $I$. Let $a>0$ be a fixed number. We assume $aI$ be the dilation
cube with the same center of $I$, and its length is $al(I)$. \vskip
0.2 cm
 \noindent {\bf Lemma 2.1} Let $-1<\alpha\leq\frac{n}{2}$. Then we
 have
quasi-norm $\|f\|_{\mathcal {L}^{2,n-2\alpha,\alpha}}$ is
well-defined. \vskip 0.2 cm
 {\bf Proof:}\ As for another bump test
function, we have the expression
$$\|f\|_{\mathcal {L}^{2,n-2\alpha,\alpha}}^{\prime}=\sup_I(\frac{1}{|I|^{1-\frac{2\alpha}{n}}}\sum_{j\geq-\log_2l(I)}2^{2\alpha
j}\|\Delta^{\prime}_jf\|_{L^2(I)}^2)^{\frac{1}{2}}$$ We let
$f=(-\triangle )^{-\frac{\alpha}{2}}g$. By the proof of Lemma 24 in
[2]. We have known that ([2], (22))
$$\frac{1}{|I|^{1-\frac{2\alpha}{n}}}\sum_{j\geq-\log_2l(I)}2^{2\alpha
j}\|(-\triangle
)^{-\frac{\alpha}{2}}\Delta^{\prime}_jg\|_{L^2(I)}^2\lesssim
\|g\|_{\mathcal {L}^{2,n-2\alpha,0}}^2$$ for any fixed cube
$I\subset\mathbb{R}^n$.

Because of proposition 8 in [2], $\mathcal
{L}^{2,n-2\alpha,0}=L^{2,n-2\alpha}$ is Campanato space and thus
well defined. We have
$$\|g\|_{L^{2,n-2\alpha}}\lesssim
\sup_I(\frac{1}{|I|^{1-\frac{2\alpha}{n}}}\sum_{j\geq-\log_2l(I)}\|\Delta_jg\|_{L^2(I)}^2)^{\frac{1}{2}}$$Then
$$\|f\|_{\mathcal
{L}^{2,n-2\alpha,\alpha}}^{\prime}\lesssim \|f\|_{\mathcal
{L}^{2,n-2\alpha,\alpha}}$$ by lemma 24 in [2]. \vskip 0.2 cm
\noindent {\bf Lemma 2.2} Let $\alpha>0$. We have another quasi-norm
definition of $\mathcal {L}^{2,n-2\alpha,\alpha}$ as
follows:$$\|f\|_{\mathcal
{L}^{2,n-2\alpha,\alpha}}=[\sup_I\sum_{k\geq0}2^{(2\alpha-n)k}\sum_{J\in
D_k(I)}\frac{1}{|J|}\sum_{j\geq-\log_2l(J)}\|\Delta_jf\|_{L^2(J)}^2]^{\frac{1}{2}}$$
\vskip 0.2 cm
 {\bf Proof:} For a fixed
$I\subset\mathbb{R}^n$,$$\sum_{k\geq0}2^{(2\alpha-n)k}\sum_{J\in
D_k(I)}\frac{1}{|J|}\sum_{j\geq-\log_2l(J)}\|\Delta_jf\|_{L^2(J)}^2$$$$=\sum_{k\geq0}\frac{1}{|I|}2^{2\alpha
k}\sum_{j\geq k-\log_2l(I)}\|\Delta_jf\|_{L^2(I)}^2$$If
$f\in\mathcal {L}^{2,n-2\alpha,\alpha}$. By Fubini theorem, we
exchange the order of summation of above identity as
follows:$$\sum_{k\geq0}\frac{1}{|I|}2^{2\alpha k}\sum_{j\geq
k-\log_2l(I)}\|\Delta_jf\|_{L^2(I)}^2=\frac{1}{|I|}\sum_{j\geq-\log_2l(I)}(\sum_{k=0}^{j+\log_2l(I)}2^{2\alpha
k})\|\Delta_jf\|_{L^2(I)}^2$$$$\ \ \ \ \ \ \ \ \ \ \ \ \ \ \ \ \ \ \
\ \ \ \ \ \ \ \ \ \ \ \ \ \ \ \ \ \ \ \ \
\thickapprox\frac{1}{|I|^{1-\frac{2\alpha}{n}}}\sum_{j\geq-\log_2l(I)}2^{2\alpha
j}\|\Delta_jf\|_{L^2(I)}^2$$Then
$$\sum_{k\geq0}2^{(2\alpha-n)k}\sum_{J\in
D_k(I)}\frac{1}{|J|}\sum_{j\geq-\log_2l(J)}\|\Delta_jf\|_{L^2(J)}^2=\frac{1}{|I|^{1-\frac{2\alpha}{n}}}\sum_{j\geq-\log_2l(I)}2^{2\alpha
j}\|\Delta_jf\|_{L^2(I)}^2<\infty\eqno(2.1)$$ On the other hand. If
$$\sup_I\sum_{k\geq0}2^{(2\alpha-n)k}\sum_{J\in
D_k(I)}\frac{1}{|J|}\sum_{j\geq-\log_2l(J)}\|\Delta_jf\|_{L^2(J)}^2<\infty$$
We have (2.1) is also valid by Fubini theorem. Then we complete the
proof.\vskip 0.2 cm

 \noindent {\bf Lemma 2.3}
Let $m\geq2$, $\alpha>-\frac{n}{2}$. We have
$$
\sum_{k\geq0}2^{(2\alpha-n)k}\sum_{J\in
D_k(I)}\frac{1}{|J|^2}\int_{mJ}\int_{mJ}|f(x)-f(y)|^2dxdy\lesssim
m^{2\alpha+2n}\|f\|_{Q_{\alpha}}^2 \eqno(2.2)$$ for any fixed cub
$I\subset\mathbb{R}^n$.
\vskip 0.2 cm

{\bf Proof:} If $m\geq2$, We also adopt the idea of lemma 5.3 in [3]
but need more complexity techniques. Observe
that$$\sum_{k\geq0}2^{(2\alpha-n)k}\sum_{J\in
D_k(I)}\frac{1}{|J|^2}\int_{mJ}\int_{mJ}|f(x)-f(y)|^2dydx$$$$=
l(I)^{2\alpha-n}\int_{mI}\int_{mI}k(x,y)|f(x)-f(y)|^2dxdy$$And we
have the following identity:
$$k(x,y)=\sum_{J\in
D(I)}\frac{\chi_{mJ}(x)\chi_{mJ}(y)}{l(J)^{2\alpha+n}}$$We let
$\mathbf{\Gamma}\triangleq\{J\in D(I):x,y\in mJ\}$. Then we get the
alternate expression of $k(x,y)$:$$k(x,y)=\sum_{J\in
\mathbf{\Gamma}}\frac{1}{l(J)^{2\alpha+n}}$$It is crucial to
estimate the magnitude of $k(x,y)$.

To begin with, we give a definition of {\bf allowed cubes}: Let $J$
be an allowed cube if there is no such dyadic subcube
$J^{\prime}\subset J$, such that $J^{\prime}\in\mathbf{\Gamma}$. We
note $\mathbf{\Gamma^a}$ be the set of allowed cubes.

We immediately conclude that all the allowed cubes disjoint each
other.

We assert$$k(x,y)=\sum_{J\in
\mathbf{\Gamma}}\frac{1}{l(J)^{2\alpha+n}}\approx\sum_{J\in
\mathbf{\Gamma^a}}\frac{1}{l(J)^{2\alpha+n}}\eqno(2.3)$$We now
prove(2.3): First, it is trivial $$\sum_{J\in
\mathbf{\Gamma}}\frac{1}{l(J)^{2\alpha+n}}\geq\sum_{J\in
\mathbf{\Gamma^a}}\frac{1}{l(J)^{2\alpha+n}}.$$For any
$J\in\mathbf{\Gamma}$, there exists only one sequence of dyadic
cubes $J_k$($k=1,...$), such that $J\subset J_1\subset
J_2\subset...$, and $J_k\in\mathbf{\Gamma}$. We define a partial
order "$<$": $J_1<J_2$ if and only if $J_1\subset J_2$. Notice that
$\mathbf{\Gamma^a}$ essentially correspond the equivalent class of
$\mathbf{\Gamma}$. We denote $\mathbf{T_{J_0}}$ be the tree which
contains $J_0$. We have the covering property:
$$\bigcup_{J\in\mathbf{\Gamma}}J\subset\bigcup_{J_0\in\mathbf{\Gamma^a}}\bigcup_{J_1\in\mathbf{T_{J_0}}}J_1$$
By $\alpha>-\frac{n}{2}$ we have following estimate:
$$\sum_{J\in
\mathbf{\Gamma}}\frac{1}{l(J)^{2\alpha+n}}\leq\sum_{J_0\in
\mathbf{\Gamma^a}}\sum_{J\in\mathbf{T_{J_0}}}\frac{1}{l(J)^{2\alpha+n}}\leq
C(n,\alpha)\sum_{J_0\in
\mathbf{\Gamma^a}}\frac{1}{l(J_0)^{2\alpha+n}}.$$This indicate (2.3)
is valid.

Having established (2.3), we turn to estimate the magnitude of
$k(x,y)$. We denote a initial cube $I_0$ with the edge parallel to
the coordinate axes and contains $x,y$. The $I_0$ is fixed and set
its length $l(I_0)=\sqrt{n}|x-y|$. Here $I_0$ does not necessary
belongs to $D(I)$. We define a sequence of cubes ${I_k}$
($k=0,1,2,...$) such that $I_k=2^kI_0$.

Then we split $\mathbf{\Gamma^a}$ into two kinds of sets. First, we
let
$$\mathbf{\Gamma_0^{(1)}}\triangleq\{J\in\mathbf{\Gamma^a}:J\cap
I_0\neq\varnothing,J\subset I_1\}.$$ When $k\geq1$, we define the
following first kind of sets
inductively:$$\mathbf{\Gamma_k^{(1)}}\triangleq\{J\in\mathbf{\Gamma^a}:J\cap
I_k\neq\varnothing,J\subset I_{k+1},
J\cap\cup_{j=0}^{k-1}I_j=\varnothing\}.$$We get first kind of sets
by induction.

The second kind of sets are the complement of the first kind of sets
counterpart. We construct these sets as follows:
Let$$\mathbf{\Gamma_0^{(2)}}\triangleq\{J\in\mathbf{\Gamma^a}:J\cap
I_0\neq\varnothing,J\nsubseteq I_1\},$$ and also
define:$$\mathbf{\Gamma_k^{(2)}}\triangleq\{J\in\mathbf{\Gamma^a}:J\cap
I_k\neq\varnothing,J\nsubseteq I_{k+1},
J\cap\cup_{j=0}^{k-1}I_j=\varnothing\}.$$The second kind of sets
then given by induction.

We can immediately deduce
$$\mathbf{\Gamma^a}=\bigcup_{k\geq0}\mathbf{\Gamma_k^{(1)}}\bigcup\mathbf{\Gamma_k^{(2)}}.$$By
(2.3),$$k(x,y)\leq\sum_{k\geq0}(\sum_{J\in
\mathbf{\Gamma_k^{(1)}}}\frac{1}{l(J)^{2\alpha+n}}+\sum_{J\in
\mathbf{\Gamma_k^{(2)}}}\frac{1}{l(J)^{2\alpha+n}})=\mathbb{I}+\mathbb{II}.\eqno(2.4)$$

The estimate of $\mathbb{I}$:

For any cube $J\in\mathbf{\Gamma_k^{(1)}}$, let $l_j
\triangleq\min\{l(J):J\in\mathbf{\Gamma_j^{(1)}}\}$, ($j\geq1$). By
geometric properties, and its definition, we know that the segment
$[x,y]$ should be contained in $mJ$. By definition of
$\mathbf{\Gamma_k^{(1)}}$, we have $\sqrt{n}l_0\geq m^{-1}|x-y|$,
and also $\sqrt{n}l_1\geq m^{-1}|x-y|$. Also, we know that $mJ$
intersects the area of $I_k\cap I_0^c$ for $k\geq2$. (See figure 1)
Then we have$$ml_{k}\geq
\frac{1}{2}(l(I_{k-1})-l(I_0))=\frac{2^{k-1}-1}{2}l(I_0)$$

Since all of the cubes in $\mathbf{\Gamma_k^{(1)}}$ contained in
$I_{k+1}$. We could calculate the number of elements in
$\mathbf{\Gamma_k^{(1)}}$:$$\#\mathbf{\Gamma_k^{(1)}}\leq\frac{l(I_{k+1})^n}{l_k^n}\leq
C_1(n)m^n$$ Thus the estimate of $\mathbb{I}$ is
clear:$$\mathbb{I}=\sum_{k\geq0}\sum_{J\in
\mathbf{\Gamma_k^{(1)}}}\frac{1}{l(J)^{2\alpha+n}}\leq
C_1(n)m^{2\alpha+2n}\sum_{k\geq0}2^{-2\alpha
k-nk}|x-y|^{-2\alpha-n}$$Because $\alpha>-\frac{n}{2}$. We could
deduce
$$\mathbb{I}\lesssim
m^{2\alpha+2n}|x-y|^{-2\alpha-n}.\eqno(2.5)$$

The estimate of $\mathbb{II}$:

For each $J\in\mathbf{\Gamma_k^{(2)}}$, notice that all of $J$
intersect the area of $I_{k+1}\cap I_k^c$. We have
$l(J)\geq\frac{1}{2}(2^{k+1}-2^k)l(I_0)$. The cross-sections $R_k$
are rectangles have the mini-length greater than $2^{k-1}l(I_0)$, or
at least contain a rectangle which has the mini-length greater than
$2^{k-1}l(I_0)$. Also, $J\in\mathbf{\Gamma_k^{(2)}}$ disjoint each
other and therefore all of $R_k$ are disjoint each other as well.
(See figure 2) We immediately obtain the number of elements in
$\mathbf{\Gamma_k^{(2)}}$
satisfies:$$\#\mathbf{\Gamma_k^{(2)}}\leq\max\{\frac{|I_{k+1}\cap
I_k^c|}{|R_k|}: R_k=J\cap I_{k+1}\cap I_k^c,
J\in\mathbf{\Gamma_k^{(2)}}\}\leq C_2(n).$$ Thus we have the
estimate of $\mathbb{II}$:$$\mathbb{II}=\sum_{k\geq0}\sum_{J\in
\mathbf{\Gamma_k^{(2)}}}\frac{1}{l(J)^{2\alpha+n}}\leq\sum_{k\geq0}C_2(n)2^{-n
k(2\alpha+n)}|x-y|^{-2\alpha-n}.$$Because $\alpha>-\frac{n}{2}$. We
have proved following
estimate:$$\mathbb{II}\lesssim|x-y|^{-2\alpha-n}.\eqno(2.6)$$
Combining estimates (2.3)(2.4)(2.5)(2.6), we get the desired
conclusion by (1.1). Notice that if there exists some $k$ or $j$
($j=0,1$) such that $\mathbf{\Gamma_k^{(j)}}=\emptyset$. It will
lead the (2.4) be a lacunary series, and this do not effect the
correctness of the results. We then complete the proof of Lemma 2.4.

\vskip 0.1 cm

\section{Proof of the main theorem}
In the following discussion, all of the cube $I\subset\mathbb{R}^n$
have the parallel to the coordinate axes edges.

\vskip 0.2 cm

\noindent {\bf The proof of statement: "If $f\in\mathcal
{L}^{2,n-2\alpha,\alpha}$ then
$\|f\|_{Q_{\alpha}}\lesssim\|f\|_{\mathcal
{L}^{2,n-2\alpha,\alpha}}$"}: \vskip 0.2 cm

For $f\in S^\prime/\mathcal {P}$, and for a fixed cube $I$, we
decompose $f$ as follows:
$$f=\sum_{j\in \mathbb{Z}}\Delta_j(f)(x)=\sum_{j<-\log_2
l(I)}\Delta_j(f)(x)+\sum_{j\geq-\log_2 l(I)}\Delta_j(f)(x)$$Then we
have$$l(I)^{2\alpha-n}\int_I\int_I\frac{|f(x)-f(y)|^2}{|x-y|^{2\alpha+n}}dx
dy$$$$\lesssim l(I)^{2\alpha-n}\int_I\int_I|\sum_{j<-\log_2
l(I)}\Delta_j(f)(x)-\sum_{j<-\log_2
l(I)}\Delta_j(f)(y)|^2|x-y|^{-2\alpha-n}dx
dy$$$$+l(I)^{2\alpha-n}\int_I\int_I|\sum_{j\geq-\log_2
l(I)}\Delta_j(f)(x)-\sum_{j\geq-\log_2
l(I)}\Delta_j(f)(y)|^2|x-y|^{-2\alpha-n}dx
dy$$$$\triangleq\mathbb{III}+\mathbb{IV}\eqno(3.1)$$

The estimate of $\mathbb{III}$:

In [2], we have known $$\sum_{j<-\log_2 l(I)}\max_{x\in
I}|\partial_x\Delta_j f(x)|\leq\|f\|_{BMO}l(I)^{-1}$$ Combining the
trivial property $\mathcal{L}^{2,n-2\alpha,\alpha}\subset BMO$ and
the fact $\alpha\in(0,1)$. We have
$$\mathbb{III}\leq\|f\|_{BMO}^2
l(I)^{2\alpha-n-2}\int_I\int_I|x-y|^{2-2\alpha-n}dxdy\lesssim\|f\|^2_{\mathcal{L}^{2,n-2\alpha,\alpha}}\eqno(3.2)$$

The estimate of $\mathbb{IV}$:

First, we rewrite$$\mathbb{IV}=l(I)^{2\alpha-n}\int_{|y|\leq
l(I)}\int_I|\sum_{j\geq-\log_2
l(I)}\Delta_j(f)(x)-\sum_{j\geq-\log_2
l(I)}\Delta_j(f)(x+y)|^2dx|y|^{-2\alpha-n}dy$$

The following arguments are rather standard as the proof of
$\|f\|_{L^{2,\lambda}}\lesssim\|f\|_{\mathcal{L}^{2,n-2\alpha,\alpha}}$
in [2], but need a slight modification.

There exists $\theta(\xi)\in C^{\infty}_0$ be a positive and radial
function such that $\check{\theta}(x)\geq1$, for $|x|\leq
\frac{1}{\pi}$ and supported in
$\{\xi\in\mathbb{R}^n:|\xi|\leq\frac{1}{2}\}$. We denote $c(I)$ be
the center of $I$. Let
$$\varphi_{I}(x)=l(I)^{\alpha-\frac{n}{2}}\check{\theta}(\pi \frac{x-c(I)}{l(I)})$$ For this fixed cube $I$, Schwartz
function $\varphi_{I}$ has the following properties:$$\ \ \ \
|\varphi_{I}(x)|^2\geq Cl(I)^{2\alpha-n},
 x\in
I$$$$\mathrm{supp}
\widehat{\varphi_{I}}(\xi)\subset\{\xi\in\mathbb{R}^n:|\xi|\leq\frac{1}{2}l(I)^{-1}\}$$
Then$$\mathbb{IV}\leq\int_{|y|\leq
l(I)}\int_{\mathbb{R}^n}|\varphi_{I}(x)|^2|\sum_{j\geq-\log_2
l(I)}\Delta_j(f)(x)-\sum_{j\geq-\log_2
l(I)}\Delta_j(f)(x+y)|^2dx|y|^{-2\alpha-n}dy\eqno(3.3)$$By
Plancherel
theorem,$$\int_{\mathbb{R}^n}|\varphi_{I}(x)|^2|\sum_{j\geq-\log_2
l(I)}\Delta_j(f)(x)-\sum_{j\geq-\log_2
l(I)}\Delta_j(f)(x+y)|^2dx$$$$=\int_{\mathbb{R}^n}|\sum_{j\geq-\log_2
l(I)}(\widehat{\varphi_{I}}(\xi)*\widehat{\Delta_j(f)}(\xi)|^2|1-e^{-2i\pi
y\xi}|^2d\xi\eqno(3.4)$$And because of $|1-e^{-2i\pi
y\xi}|\leq\min\{2,C_{\mu_0}|y|^{\mu_0}|\xi|^{\mu_0}\}$. We note
$\mu_0$ be a fixed positive number with $\alpha<\mu_0<1$. We have
the fact$$\int_{|y|\leq l(I)}\frac{|1-e^{-2i\pi
y\xi}|^2}{|y|^{2\alpha+n}}dy\lesssim|\xi|^{2\mu_0}\int_{|y||\xi|\leq1}|y|^{2\mu_0-2\alpha-n}dy+\int_{|y||\xi|\geq1}|y|^{-2\alpha-n}dy\lesssim
|\xi|^{2\alpha}\eqno(3.5)$$We define another Littlewood-Paley
operator:$$\widehat{\Delta_j^{\prime}(f)}(\xi)=|2^{-j}\xi|^{\alpha}\widehat{\psi}(2^{-j}\xi)\widehat{f}(\xi)\eqno(3.6)$$Because
of the orthogonality property, we citing the following estimate in
[2]
$$|\sum_{j\geq-\log_2
l(I)}(\widehat{\varphi_{I}}(\xi)*\widehat{\Delta_j^{\prime}(f)}(\xi)|^2\leq7\sum_{j\geq-\log_2
l(I)}|(\widehat{\varphi_{I}}(\xi)*\widehat{\Delta_j^{\prime}(f)}(\xi)|^2\eqno(3.7)$$Combining
(3.3)(3.4)(3.5)(3.6)(3.7) as well as exchange the order of
integration of (3.3), we have
$$\mathbb{IV}\lesssim\sum_{j\geq-\log_2
l(I)}\int_{\mathbb{R}^n}|\widehat{\varphi_{I}}(\xi)*\widehat{\Delta_j^{\prime}(f)}(\xi)|^22^{2\alpha
j} d\xi=\sum_{j\geq-\log_2
l(I)}\int_{\mathbb{R}^n}|\varphi_{I}(x)\Delta_j^{\prime}(f)(x)|^22^{2\alpha
j}dx.$$ The following arguments are almost the same as in [2].

Denote $k\in \mathbb{Z}^n$,
$a_k=\max\{|\check{\theta}(x)|^2:|x-k|\leq\frac{1}{2}\}$. We let
$Q_k$ be the disjoint cubes in $\mathbb{R}^n$ have the center at
$l(I)k$ with the length of $l(I)$. Then $Q_k$ ($k\in\mathbb{Z}^n$)
become the partition of $\mathbb{R}^n$. We have
$$\mathbb{IV}\lesssim \sum_{k\in \mathbb{Z}^n}a_k\sum_{j\geq-\log_2
l(I)}\frac{1}{|I|^{1-\frac{2\alpha}{n}}}\int_{Q_k}|\Delta_j^{\prime}(f)(x)|^22^{2\alpha
j}dx$$By the property of Schwartz function and Lemma 2.1 we
have$$\mathbb{IV}\lesssim\sup_I\frac{1}{|I|^{1-\frac{2\alpha}{n}}}\sum_{j\geq-\log_2l(I)}2^{2\alpha
j}\|\Delta^{\prime}_jf\|_{L^2(I)}^2\lesssim\|f\|_{\mathcal
{L}^{2,n-2\alpha,\alpha}}^2$$Combining above estimate and
(3.1)(3.2), we have
$$l(I)^{2\alpha-n}\int_I\int_I\frac{|f(x)-f(y)|^2}{|x-y|^{2\alpha+n}}dx
dy\lesssim\|f\|_{\mathcal {L}^{2,n-2\alpha,\alpha}}^2$$for any fixed
cube $I$.

By (1.1) we complete the proof of
$\|f\|_{Q_{\alpha}}\lesssim\|f\|_{\mathcal
{L}^{2,n-2\alpha,\alpha}}$.

\vskip 0.2 cm

\noindent {\bf The proof of statement: "If $f\in
Q_{\alpha}(\mathbb{R}^n)$, then $\|f\|_{\mathcal
{L}^{2,n-2\alpha,\alpha}}\lesssim\|f\|_{Q_{\alpha}}$" }:

\vskip 0.2 cm

To begin with, by lemma 2.2, it suffices to show
$$\sum_{k\geq0}2^{(2\alpha-n)k}\sum_{J\in
D_k(I)}\frac{1}{|J|}\sum_{j\geq-\log_2l(J)}\|\Delta_jf\|_{L^2(J)}^2\lesssim\|f\|_{Q_{\alpha}}^2\eqno(3.8)$$for
any fixed cube $I$.

For any fixed subcube $J\subset I$, we have the decomposition of $f$
related to $J$ as follows:
$f=(f-f_{2J})\chi_{2J}+(f-f_{2J})\chi_{(2J)^c}+f_{2J}$. Then we have
the following decomposition:
$$\sum_{k\geq0}2^{(2\alpha-n)k}\sum_{J\in
D_k(I)}\frac{1}{|J|}\sum_{j\geq-\log_2l(J)}\|\Delta_jf\|_{L^2(J)}^2$$$$\lesssim\sum_{k\geq0}2^{(2\alpha-n)k}\sum_{J\in
D_k(I)}\frac{1}{|J|}\sum_{j\geq-\log_2l(J)}\|\Delta_j(f-f_{2J})\chi_{2J}\|_{L^2(J)}^2$$$$+\sum_{k\geq0}2^{(2\alpha-n)k}\sum_{J\in
D_k(I)}\frac{1}{|J|}\sum_{j\geq-\log_2l(J)}\|\Delta_j(f-f_{2J})\chi_{(2J)^c}\|_{L^2(J)}^2$$$$+\sum_{k\geq0}2^{(2\alpha-n)k}\sum_{J\in
D_k(I)}\frac{1}{|J|}\sum_{j\geq-\log_2l(J)}\|\Delta_jf_{2J}\|_{L^2(J)}^2\triangleq\mathbb{V}+\mathbb{VI}+\mathbb{VII}$$
It is obviously that $f_{2J}$ is a constant and we have
$\Delta_j(f_{2J})\equiv0$ for all the $j\in\mathbb{Z}$ and all the
subcube $J\subset I$. Then we have $\mathbb{VII}=0$. In order to
prove (3.8), we only need to demonstrate
$\mathbb{V}\lesssim\|f\|_{Q_{\alpha}}^2$ and also
$\mathbb{VI}\lesssim\|f\|_{Q_{\alpha}}^2$.

The estimate of $\mathbb{V}$:

By Plancherel theorem, we have
$$\sum_{j\in\mathbb{Z}}\int_{R^n}|\Delta_j(f-f_{2J})\chi_{2J}|^2dx=\|
(f-f_{2J})\chi_{2J}\|_{L^2}^2.$$ We can
deduce$$\mathbb{V}\leq\sum_{k\geq0}2^{(2\alpha-n)k}\sum_{J\in
D_k(I)}{\frac{1}{|J|}}\sum_{j>-\log_2
l(J)}\int_{\mathbb{R}^n}|\Delta_j(f-f_{2J})\chi_{2J}|^2dx$$$$=\sum_{k=0}^{\infty}2^{(2\alpha
-n)k}\sum_{J\in D_k(I)}\frac{1}{|J|}\int_{2J}| f-f_{2J}| ^2dx$$Then
$\mathbb{V}\lesssim\|f\|_{Q_\alpha}^2$ follows by (2.2) with the
case of $m=2$.

The estimate of $\mathbb{VI}$:

To start with, we assume $x\in J$. We give the following arguments:
$$|\Delta_j(f-f_{2J})\chi_{(2J)^c}(x)|=|\int_{R^n}\psi_j(x-y)(f(y)-f_{2J})\chi_{(2J)^c}(y)dy|$$
$$\leq
\sum_{l\geq1}\int_{2^{l+1}J\bigcap(2^lJ)^c}|\psi_j(x-y)||f(y)-f_{2J}|dy\eqno(3.9)$$Since
$\psi$ is a Schwartz function, then $\psi$ descend faster than any
polynomial. Let $M > 2\alpha + n$ be a fixed large number. We have
$$|\psi_j(x-y)|\leq C_M2^{jn}(1+|x-y|2^j)^{-M-n}\eqno(3.10)$$Notice
that $|x-y|\geq2^{l-1}l(J)$. By (3.9)(3.10), we have the
Littlewood-Paley operator could be controlled by the mean
oscillation:$$|\Delta_j(f-f_{2J})\chi_{(2J)^c}(x)|\lesssim
2^{-jM}l(J)^{-M}\sum_{l\geq1}2^{-lM}(|f-f_{2J}|)_{2^{(l+1)}J}\eqno(3.11)$$We
could also deduce the following estimate by Cauchy-Schwarz
inequality
$$\sum_{l\geq1}2^{-lM}(|f-f_{2J}|)_{2^{(l+1)}J}\leq (\sum_{l\geq1}2^{-lM})^{\frac{1}{2}}(\sum_{l\geq1}2^{-lM}(|f-f_{2J}|)_{2^{(l+1)}J}^2)^{\frac{1}{2}}\eqno(3.12)$$
Combining (3.11)(3.12) and using the Jensen inequality, we get the
estimate of $\mathbb{VI}$ as follows:$$\mathbb{VI}\lesssim
\sum_{k\geq0}2^{(2\alpha-n)k}\sum_{J\in
D_k(I)}{\sum_{l\geq1}2^{-lM}(\frac{1}{|2^{l+1}J|}\int_{2^{l+1}J}|f-f_{2J}|dy})^2$$$$\
\ \ \ \ \lesssim\sum_{l\geq1}2^{-lM-ln}\sum_{k\geq0}\sum_{J\in
D_k(I)}\frac{1}{|J|^2}\int_{2^{l+1}J}\int_{2^{l+1}J}|f(x)-f(y)|^2dxdy$$Using
the growth estimate provided in Lemma 2.3. The above summation could
be exchanged and we could obtain
$$\mathbb{VI}\lesssim\sum_{l\geq1}2^{-l(M-2\alpha-n)}\|f\|_{Q_{\alpha}}^2\lesssim\|f\|_{Q_{\alpha}}^2$$This
completes the proof.

\section{Remark}

In fact, we have known that $Q_{\alpha}(\mathbb{R}^n)\subset\mathcal
{L}^{2,n-2\alpha,\alpha}$ for $-\infty<\alpha<\infty$. But $\mathcal
{L}^{2,n-2\alpha,\alpha}\subset Q_{\alpha}(\mathbb{R}^n)$ probably
no longer available for $\alpha\geq1$. That means if $\alpha\geq1$,
$f\in{L}^{2,n-2\alpha,\alpha}$. Then we cannot deduce $f(x)$ is a
constant function. At least, if we let $\alpha=1$, $n=2$. We could
easily construct a non-constant Sobolev function $f(x)$, such that
$\partial_x f(x)\in L^2(R^2)$. For example, let $f(x)$ be a
non-constant Schwartz function. By ([2], theorem 10), we know that
$f\in\mathcal {L}^{2,0,1}$.


\begin{thebibliography}{99}

\bibitem{l1} R. Aulaskari, J. Xiao, R. Zhao,  On subspaces and subsets of BMOA and UBC. Analysis,
15 (1995), 101-121.

\bibitem{i1} A. El. Baraka, Littlewood-Paley characterization for Campanato spaces. J. Function Spaces and Applications 4,
No.2 (2006), 193-220.

\bibitem{i1} G. Dafni, J. Xiao Some new tent spaces and duality theorems for
fractional Carleson measures and $Q_\alpha(\Bbb R^n)$.  J. Funct.
Anal., 208 (2004), 377-422.

\bibitem{c1} M. Ess\'{e}n, S. Janson, Lizhong. Peng, Jie. Xiao,  $Q$ Spaces of several real variables. Indiana Univ. Math. J. 49, No.2 (2000),
575-615.

\bibitem{m1}  M. Ess\'{e}n, J. Xiao,
 Some results on $Q_p$ spaces, $0<p<1$. J. Reine Angew. Math.
485 (1997), 173-195.

\bibitem{r1} H. Kozono and M. Yamazaki, Semilinear heat equations and the Navier-Stokes equations
with distributions in new function spaces as initial data. Comm.
Partial Differential Equations, 19 (1994), 959-1014.

\bibitem{v2} J. Xiao, Homothetic variant of fractional Sobolev space with application to Navier-Stokes system.  Dynamics of P.D.E. 4 (2007), 227-245.


\end{thebibliography}
\end{document}